\documentclass[11pt,a4paper,twoside]{amsart}
\usepackage{geometry}
\usepackage{amsfonts, amssymb, amsmath, amscd, latexsym}
\usepackage{psfrag, eucal, enumerate, latexsym}
\usepackage{amsthm}
\usepackage{xypic}
\usepackage{enumerate}
\usepackage{psfrag}
\usepackage{graphicx}
\theoremstyle{proposition}
\newtheorem*{teo}{Theorem}
\newtheorem{theorem}{Theorem}[section]
\newtheorem{lemma}[theorem]{Lemma}
\newtheorem{definition}[theorem]{Definition}
\theoremstyle{remark}
\newtheorem{remark}[theorem]{Remark}
\numberwithin{equation}{section}
\theoremstyle{corollary}
\newtheorem{corollary}[theorem]{Corollary}
\theoremstyle{proposition}
\newtheorem{proposition}[theorem]{Proposition}

\theoremstyle{remark}
\newtheorem{example}[theorem]{Example}
\newcommand{\CC}{\mathbb C}
\newcommand{\NN}{\mathbb N}
\newcommand{\PP}{{\mathbb P}}

\newcommand{\ZZ}{{\mathbb Z}}

\newcommand{\OO}{\mathcal{O}}

\DeclareMathOperator{\HH}{H}
\DeclareMathOperator{\hh}{h}
\DeclareMathOperator{\Ext}{Ext}
\DeclareMathOperator{\Hom}{Hom}
\DeclareMathOperator{\im}{Im}
\DeclareMathOperator{\rk}{rk}
\DeclareMathOperator{\Ker}{Ker}

\DeclareMathOperator{\coker}{Coker}

\DeclareMathOperator{\I}{Id}
\DeclareMathOperator{\Stab}{Stab}

\DeclareMathOperator{\Ad}{Ad}

\DeclareMathOperator{\GL}{GL}

\DeclareMathOperator{\SL}{SL}
\DeclareMathOperator{\End}{End}

\DeclareMathOperator{\Sym}{Sym}
\begin{document}
%
\title{Cokernel bundles and Fibonacci bundles}

\author{Maria Chiara Brambilla}
\address{Dipartimento di Matematica
e Applicazioni per l'Architettura,
Universit\`a di Firenze,
piazza Ghiberti, 27 - 50122, Florence, Italy}
\thanks{Corresponding author: e-mail: {\sf brambilla@math.unifi.it} 
}

\begin{abstract}
{We are interested in those bundles $C$ on $\PP^N$ which admit a resolution
of the form
\[0\rightarrow \CC^s\otimes E\xrightarrow{\mu}
\CC^t\otimes F\rightarrow C\rightarrow 0.\] In this paper we prove
that, under suitable conditions on $(E,F)$, a generic bundle with
this form is either simple or canonically decomposable. As
applications we provide an easy criterion for the stability of
such bundles on $\PP^2$ and we prove the stability when $E=\OO,
F=\OO(1)$ and $C$ is an exceptional bundle on $\PP^N$ for
$N\geq2$.}
\end{abstract}
%
%
\maketitle
\section{Introduction}
In this paper we consider the family of those vector bundles $C$
on $\PP^N$, with $N\geq2$, which admit a resolution of the form
\begin{equation}\label{SD}
0\rightarrow E^s\xrightarrow{\mu} F^t\rightarrow C\rightarrow 0,
\end{equation}
for some bundles $E,F$ and for $s,t\in\NN$.
Here and in the sequel we write $E^s$ (resp.\ $F^t$) instead of
$\CC^s\otimes E$ (resp.\ $\CC^t\otimes F$), and we assume  $t\rk(F)-s\rk(E)\geq N$.
Any bundle $C$ in \eqref{SD} is the cokernel of a morphism of bundles $\mu$.
Our purpose is to describe the properties of the bundles corresponding to
generic morphisms in $\Hom(E^s, F^t)$.
In particular we want to find criteria of simplicity, rigidity and decomposability.

Throughout the paper we will assume that $E$ and $F$ are two different
vector bundles on $\PP^N=\PP(V)$, with $N\geq2$, which verify the following
{\em basic hypotheses}:
\begin{eqnarray}
&E \mbox{ and }F \mbox{ are simple, and }\Hom(F,E)=\Ext^1(F,E)=0,\label{ipo1}\\
&\begin{small}
\mbox{ the sheaf }E^*\otimes F\mbox{ is globally generated,
and }W=\Hom(E,F) \mbox{ has dimension }w\geq3.
\end{small}\label{ipo2}
\end{eqnarray}

The first instance one can consider is $E=\OO$ and $F=\OO(1)$: in this case
we obtain the family of {\em Steiner bundles}, where the morphism $\mu$ is
a $(t\times s)$-matrix whose entries are homogeneous linear polynomials.
In \cite{mio} we studied this case and we obtained a criterion for the
simplicity of Steiner bundles. In \cite{tesi} we extended this result,
describing the canonical decomposition of generic non-simple Steiner
bundles. In particular we proved that the indecomposable elements which appear in
such a decomposition are exceptional bundles.

{\em Exceptional bundles} were introduced by Dr\'ezet and Le Potier in
\cite{DrezetLePotier} as a class of bundles on $\PP^2$ without
deformation.
Later the school of Rudakov generalized the concept of exceptional bundles and
introduced mutations in order to construct them, in the setting of derived categories (see for
example \cite{Rudakovseminari}).

The first result we obtain here (Theorem \ref{simpli} below) states that if
$C$ is the cokernel of a
generic map $\mu\in\Hom(E^s, F^t)$, then
\[C\ \textrm{is simple}\quad \Leftrightarrow\quad s^2-wst+t^2\leq1.\]

This result allows us to obtain a criterion for the stability of the
cokernel bundles $C$ on $\PP^2$.
In fact Dr\'ezet and Le Potier obtained an important criterion
for the stability of all bundles on $\PP^2$ (see
\cite{DrezetLePotier}), but their result is very difficult to
apply. In this paper, using another result of Dr\'ezet (see
\cite{Drezet99}),
we get a new criterion for the stability of the
bundles $C$ with resolution \eqref{SD} on $\PP^2$, which is much easier to apply.

Our second result is a  canonical decomposition for non-simple cokernel bundles.
In this context, the main tool is a new family of bundles, here referred to as
{\em Fibonacci bundles}, which play the role of the exceptional bundles,
but which are much more general.
In Theorem \ref{mutazioni} we define Fibonacci bundles by means of mutations,
and in Theorem \ref{risoll} we prove that they admit a resolution \eqref{SD}
in which the coefficients are related to the numbers of Fibonacci
(this motivates our choice of their name).

Under additional conditions on $(E,F)$ we prove that all the Fibonacci
bundles are simple and rigid.
These two crucial properties allow us to find a canonical
decomposition of a generic non-simple cokernel bundle $C$ in
\eqref{SD} in terms of Fibonacci bundles.
More precisely, we add the following conditions on the pair $(E,F)$:
\begin{equation}\label{ipo3}
E\mbox{ and }F\mbox{ are rigid, and }\Ext^2(F,E)=0,
\end{equation}
and we call, for brevity, {\em hypotheses $(R)$} the union of conditions \eqref{ipo1},
\eqref{ipo2} and \eqref{ipo3}.
Theorem \ref{decompo} states that
if $(E,F)$ satisfies $(R)$,
then for a generic $C$ in \eqref{SD} we have
\[s^2-wst+t^2\geq1\quad\Rightarrow\quad C\cong C_k^n\oplus C_{k+1}^m,\]
where $C_k,C_{k+1}$ are Fibonacci bundles and $n,m\in\NN$.
We stress that, in this case, any generic non-simple
bundle is rigid and homogeneous.

Finally, as an application of our results, we prove the following
\begin{teo}
Any exceptional Steiner bundle on $\PP^N$ is stable for all $N\geq2$.
\end{teo}
We recall that exceptional bundles are known to be stable on $\PP^2$
(\cite{DrezetLePotier}) and on $\PP^3$ (\cite{Zube}),
but the stability of exceptional bundles on $\PP^N$ with
$N>3$ is an open problem.

The plan of the article is as follows:
in Section $2$ we  present some basic examples, and
in Section $3$ we recall the case of Steiner bundles and their
interpretation in terms of matrices.
Section $4$ is devoted to the criterion for simplicity and
Section $5$ to Fibonacci bundles.
In Section $6$ we give the decomposition theorem for non-simple
bundles and
in Section $7$ we describe some applications of our results.
Finally Section $8$ is devoted to our results on stability.

\section{Preliminaries}
For a fixed $N\geq 2$, we are interested in the vector bundles $C$ on $\PP^N=\PP(V)$ which admit a
resolution of the form \eqref{SD}
for some bundles $E,F$ which satisfy the basic
hypotheses \eqref{ipo1} and \eqref{ipo2}
and $s,t\in \NN$
such that $t\rk(F)-s\rk(E)\geq N$.

We say that $C$ is {\em generic} when the morphism
$\mu$ is generic in the space \linebreak
$H=\Hom(E^s,F^t)\cong\CC^s\otimes \CC^t\otimes\Hom(E,F)$.
The morphism $\mu$ can be represented by a $(t\times s)$-matrix, whose
entries are morphisms from $E$ to $F$.

Let us see some examples. As in \eqref{ipo2} we denote by $W$ the vector space $\Hom(E,F)$ and by $w$ its
dimension.
\begin{example}
If $E=\OO$ and $F=\OO(d)$, it is easy to check that conditions \eqref{ipo1}
and \eqref{ipo2} are satisfied for any $d\geq1$.
Hence we deal with bundles with resolution
\begin{equation}\label{primoesempio}
0\rightarrow \OO^s\xrightarrow{\mu} \OO(d)^t\rightarrow C\rightarrow
0,
\end{equation}
where $\mu$ is a matrix whose entries are homogeneous
polynomials of degree $d$. In this case $W=\HH^0(\OO(d))=S^dV$
and $w=\binom{N+d}{d}$.
In particular when $d=1$ we obtain the case of Steiner bundles, studied
in \cite{mio}.
\end{example}

\begin{example}
For any $p\geq0$, let us denote $\Omega^p(p)=\wedge^p\Omega^1(1)$.
Given $0\leq p<N$, we consider
$E=\OO(-1)$ and $F=\Omega^p(p)$ and we obtain bundles
of the form
\begin{equation}\label{secondoes}
0\rightarrow \OO(-1)^s\xrightarrow{\mu} \Omega^p(p)^t\rightarrow C\rightarrow
0.
\end{equation}
In this case $W=\wedge^{N-p}V^{*}$, $w=\binom{N+1}{N-p}$ and the entries of the matrix $\mu$
are $(N-p)$-forms.
Analogously we consider
$E=\Omega^p(p)$ and $F=\OO$, where $0<p\leq N$, and we obtain bundles
of the form
\begin{equation}\label{secondoes-bis}
0\rightarrow \Omega^p(p)^s\xrightarrow{\mu} \OO^t\rightarrow C\rightarrow
0,\end{equation}
where $\mu$ is a matrix of $p$-forms.
\end{example}

\begin{example}
On $\PP^2=\PP(V)$ we denote  $Q=T(-1)=\Omega^1(2)$, i.e.\
\[0\rightarrow \OO(-1)\rightarrow \OO\otimes V^{*}\rightarrow
Q\rightarrow 0,\]
and $S^pQ(d)=\Sym^pQ\otimes\OO(d)$, where $p\geq1$ and $d\in\ZZ$.
Let $E=S^pQ$ and $F=S^rQ(d)$, for some fixed $p,r\geq1$ and $d\in\ZZ$, and consider
the bundles $C$ of the form
\begin{equation}\label{terzoesempio}
0\rightarrow (S^pQ)^s\xrightarrow{\mu} (S^rQ(d))^t\rightarrow C\rightarrow
0.\end{equation}
In this case $E^*\otimes F= S^p(Q^*)\otimes S^rQ(d)\cong S^pQ\otimes
S^rQ\otimes\OO(d-p)$, and hypotheses \eqref{ipo1} and \eqref{ipo2} hold true if $d> p+1$.
\end{example}

\subsection{The Fibonacci sequence}
Given any integer $w\geq3$, we introduce the following  sequence of numbers:
\begin{equation}
\label{nufibo}
a_{w,k}=\frac{ \left(\frac{w+\sqrt{w^2-4}}{2}\right)^k-
\left(\frac{w-\sqrt{w^2-4}}{2}\right)^k}{\sqrt{w^2-4}},
\end{equation}
for $k\geq 0$. This sequence satisfies the recurrence
\[\left\{
\begin{array}{l}
  a_{w,0}=0, \\
  a_{w,1}=1,\\
  a_{w,k+1}=w a_{w,k}-a_{w,k-1}, \\
\end{array}
\right.\]
In the following for brevity we will write $a_k$ instead of $a_{w,k}$,
when the value of $w$ is clear from the context.
\begin{remark}\label{rema}
In the case $w=3$, the sequence $\{a_{w,k}\}$ is exactly the odd part
of the well known Fibonacci sequence. Also if $w>3$ the numbers
$a_{w,k}$ satisfy some good relations, analogously to Fibonacci
numbers. More precisely, for any fixed $w\geq3$, we can easily prove
by induction that the following equalities hold for all $k\geq 1$:
\begin{eqnarray}
&&a_{k-1}^2+a_{ k}^2-wa_{ k-1}a_{k}=1\label{formuladue}\\
&&a_{k}^2-a_{ k+1}a_{ k-1}=1 \label{formulauno}\\
&&a_{k+1}a_{k}-a_{k-1}a_{k+2}=w\label{formulatre}
\end{eqnarray}
From \eqref{formulauno}, it also follows that $(a_{k},a_{k-1})=1$, for
all $k\geq1$.
\end{remark}
\begin{remark}\label{pellfermat}
It is possible to prove that the pairs $(s,t)= (a_{k},a_{k+1})$ are
the unique integer solutions of the diophantine equation $s^2+t^2-w
st=1$. For more details see Lemma $3.4$ of \cite{mio}.
\end{remark}

\subsection{Exceptional bundles}
Exceptional bundles were defined by Dr\'ezet and Le Potier in
\cite{DrezetLePotier} as a class of bundles on $\PP^2$ without
deformation.
These bundles appeared as some sort of exceptional points in
the study of the stability of bundles on $\PP^2$. Dr\'ezet and Le
Potier showed that these vector bundles 
are uniquely determined by
their slopes, and they described the set of all the possible slopes.
Later, the school of Rudakov (see for example
\cite{Rudakovseminari}) generalized the definition of exceptional
bundles on $\PP^N$ and other varieties, with an axiomatic
presentation in the setting of derived categories.
Following Gorodentsev and Rudakov  (\cite{GoroRudakov}) we give the
following definition:
\begin{definition}
A bundle $E$ on $\PP^N$ is {\em exceptional} if
\[\Hom(E,E)=\CC\qquad\textrm{and}\qquad \Ext^i(E,E)=0\quad \textrm{ for all } i>1.\]
\end{definition}
We recall that a bundle is called {\em semi-exceptional} when it is a direct sum of
exceptional bundles.
\section{Steiner bundles and matrices}
In this section we recall some results concerning Steiner bundles on
$\PP^N=\PP(V)$, with $N\geq2$, i.e.\ the bundles $S$ which admit a resolution of the
form
\begin{equation}\label{stein}
0\rightarrow \OO^s\xrightarrow{\mu} \OO(1)^t\rightarrow S\rightarrow
0,
\end{equation}
for some $t-s\geq N$. In this case $\mu$ belongs to the space
$H=\CC^s\otimes \CC^t\otimes V,$ which can be seen as the
space of $(t\times s)$-matrices whose
entries are homogeneous linear forms in $N+1$ variables or,
alternatively, as the space of $(s\times t\times (N+1))$-matrices
of numbers.

We consider the following action of $\GL(s)\times\GL(t)$ on $H$:
\[\GL(s)\times\GL(t)\times H\rightarrow H\]
\[(A,B,\mu)\mapsto B^{-1} \mu A.\]
Given $\mu\in H$, we denote by $(\GL(s)\times\GL(t))\mu$ the orbit of $\mu$ and
by $\Stab(\mu)$ the stabilizer of $\mu$ with respect to the action of $\GL(s)\times\GL(t)$.

In order to describe the orbits of this action, we introduce the
following definitions concerning multidimensional matrices.

We say that two matrices $\mu,\mu'\in H$ are
$\GL(s)\times\GL(t)$-equivalent if they are in the same orbit with
respect to the action of $\GL(s) \times \GL(t)$ on $H$.
This corresponds to perform Gaussian elimination on a $(s\times
t)$-matrix of linear polynomials.
\begin{definition}\label{cano}
If  $a_k=a_{N+1,k}$ is the sequence defined in
\eqref{nufibo},
we call {\em block of type $B_k$}
a matrix in
$\CC^{a_{k-1}}\otimes\CC^{a_k}\otimes \CC^{N+1}$.
Given $n,m\in\NN$,
let $s=na_{k-1}+ma_{k}$ and
$t=n{a_{k}}+m{a_{k+1}}$.
We say that a matrix $\mu\in \CC^s\otimes\CC^t\otimes\CC^{N+1}$ is
a {\em canonical matrix} if there exist decompositions
\[\CC^s=\CC^{na_{k-1}}\oplus\CC^{ma_k}\quad \textrm{ and }\quad
\CC^t=\CC^{na_{k}}\oplus\CC^{ma_{k+1}},\]
such that the matrix $\mu$ is zero except for $n$ blocks of type $B_k$
and $m$ blocks of type $B_{k+1}$ on the diagonal. We denote such a
matrix by $B_k^n\oplus B_{k+1}^m$.
\end{definition}

The following theorem describes the elements of $H$ with respect to the action above.
For the proof we refer to \cite{mio}, \cite{tesi}, and to Theorem $4$ of \cite{kac}.
\begin{theorem}\label{megakac}
Let  $H= \CC^s\otimes\CC^t\otimes \CC^{N+1}$ be endowed with the natural action of
$\GL(s)\times\GL(t)$.

$\bullet$
If $s^2+t^2-(N+1)st\leq1$, then the stabilizer of a generic element of $H$
has dimension $1$.
In particular if $s^2+t^2-(N+1)st=1$, there is a dense orbit in $H$.

$\bullet$ If $s^2+t^2-(N+1)st\geq 1$, a generic element of $H$ is
$\GL(s)\times\GL(t)$-equivalent to a canonical matrix $B_k^n\oplus B_{k+1}^m$
for unique $n,m,k\in\NN$.
\end{theorem}

\begin{remark}
After \cite{mio} and \cite{tesi} have been written, we learned that our
results on matrices turn out to be connected to a theorem of
Kac, framed in the setting of quiver theory.
More precisely, in \cite{kac} the quiver with two vertices and $w$ arrows
from the first vertex to the second one is considered, and a representation of this
quiver is exactly a $w$-uple of linear maps from one vector space into
another.
In Theorem $4$ of \cite{kac}, Kac describes the isomorphism classes of representations of
this quiver. Notice that the proofs given in \cite{mio} and \cite{tesi}
are independent from techniques of quiver theory.
\end{remark}

\psfrag{s}{$s$}
\psfrag{t}{$t$}
\psfrag{n}{$N+1$}
\begin{figure}
\begin{center}
\includegraphics[scale=0.7]{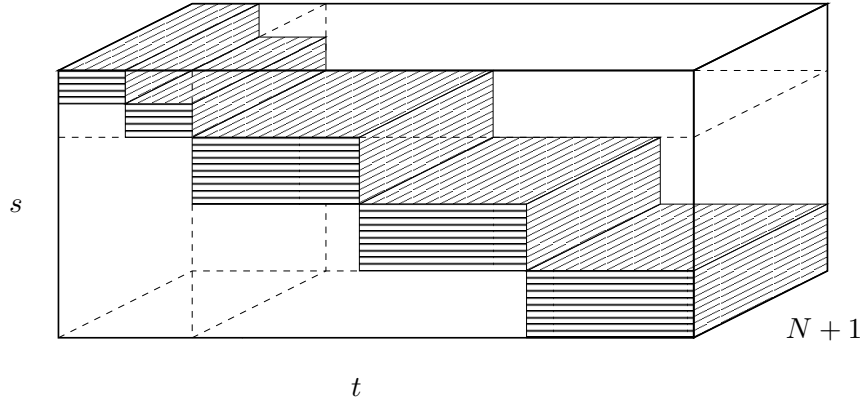}
\caption{Example of canonical matrix, with $n=2, m=3$}
\end{center}
\end{figure}
The previous theorem implies the following classification of
Steiner bundles, proved in \cite{mio} and \cite{tesi}. Here we
omit the proof, since we will prove the same result in a more general
framework later (see Theorem \ref{simpli} and Theorem \ref{decompo}).
Given $h\in\NN$ we denote by
$a_h=a_{N+1,h}$ and by $S_h$ the exceptional Steiner bundle with resolution
\[0\rightarrow \OO^{a_{h-1}}\rightarrow \OO(1)^{a_h}\rightarrow S_h\rightarrow
0.\]
\begin{theorem}\label{steiner}
Let $S$ be a generic Steiner bundle on $\PP^N$ with resolution
\eqref{stein} and  $t-s\geq N\geq2$.

$\bullet$ If $s^2-(N+1)st+t^2\leq 0$ i.e.\ if $t\leq
(\frac{N+1+\sqrt{(N+1)^2-4}}{2})s$, then the bundle $S$ is simple,

$\bullet$
if $s^2-(N+1)st+t^2\geq 1$ i.e.\ if $t> (\frac{N+1+\sqrt{(N+1)^2-4}}{2})s$, then
the bundle $S$ is isomorphic to $S_k^n\oplus S_{k+1}^m$, for
some unique $n,m,k\in\ZZ$.
\end{theorem}
Notice that the bundles of the form $S_k^n\oplus S_{k+1}^m$, which
appear in the previous theorem, correspond to canonical matrices.

%
%
\section{Simplicity}
In this section we study the simplicity of the cokernel
bundles $C$ with resolution \eqref{SD} on $\PP^N$.
As in Section $3$, we consider the natural action of $\GL(s)\times\GL(t)$ on the space $H=\Hom(E,F)$
and we denote by $\Stab(\mu)$ the stabilizer of $\mu$.

\begin{lemma}\label{lemmaStab}
If $C$ is a bundle with resolution \eqref{SD}
and $\dim \Stab(\mu)=1$, then $C$ is simple.
\end{lemma}
\begin{proof}
Assume by contradiction that $C$ is not simple. Then there exists
$\phi:C\rightarrow C$ non-trivial. Applying the functor
$\Hom(-,C)$ to the sequence \eqref{SD},
we get that $\phi$ induces $\widetilde{\phi}$ non-trivial in
$\Hom(F^t,C)$.

Now applying the functor $\Hom(F^t,-)$ again to the same sequence and using
hypothesis \eqref{ipo1}, we get 
$\Hom(F^t,F^t)\cong \Hom(F^t,C)$, hence
$\widetilde{\phi}$ induces a non-trivial morphism in $\Hom(F^t,F^t)$.
Since $F$ is simple, this non-trivial morphism induces a complex $(t\times t)$-matrix
$B\neq\lambda\I$, such that the following diagram commutes:
\[
\xymatrix{
0\ar[r]&E^s\ar[r]^{\mu}&F^t\ar[r]\ar[d]^{B}\ar[dr]^-{\widetilde{\phi}}&C\ar[r]\ar[d]^{\phi}&0 \\
0\ar[r]&E^s\ar[r]^{\mu}&F^t\ar[r]&C\ar[r]&0
}
\]
By restricting $B$ to $E^s$ and by the simplicity of $E$ we obtain a complex $(s\times s)$-matrix $A$,
such that $BM=MA$.
Let $0\neq\rho\in\CC$ be different from any eigenvalue of $B$ and $A$. If we define
$\widetilde A=A-\rho\I$ and $\widetilde B=B-\rho \I$, we get that
the pair $(\widetilde A,\widetilde B)$ belongs to $\Stab(M)\subset\GL(s)\times\GL(t)$.
Since $\widetilde A$ is not a scalar matrix, it follows that $(\widetilde A,\widetilde B)\neq(\lambda\I,\lambda\I)$
i.e.\ $\dim\Stab(M)>1$, which is a contradiction.
\end{proof}

\begin{lemma}\label{dimensioni}
If $C$ is a bundle with resolution \eqref{SD} and
$\Hom(C, F)=0$, then \linebreak
$\dim \Stab(\mu)=\dim\Hom(C,C).$
\end{lemma}
\begin{proof}
Let us apply the functor $\Hom(-,E^s)$ to the sequence \eqref{SD}.
By hypothesis \eqref{ipo1}, we obtain the following relation
\[\End\CC^s= \CC^s\otimes {\CC^s}^{*}\otimes\Hom(E,E )=
\Hom(E^s,E^s)= \Ext^1(C, E^s).\]

By applying $\Hom(-, F^t)$ to \eqref{SD} and using
hypothesis  $\Hom(C,F)=0$ and the simplicity of $F$, we get
\[0\rightarrow \End \CC^t \rightarrow\CC^t\otimes
\CC^s\otimes\Hom(E, F)\rightarrow
\CC^t\otimes\Ext^1(C,F),\]
and applying $\Hom(C,-)$ to \eqref{SD} we get
\[0\rightarrow \Hom(C,C)\rightarrow \Ext^1(C,E^s)\rightarrow \Ext^1(C,F^t).\]

The previous results together give the following commutative diagram
\[
\xymatrix@C-1ex{
&&& 0 \ar[d] &&\\
&&&\End \CC^t\ar[d]^{r_\mu}&&  \\
&&&\CC^t\otimes {\CC^s}\otimes \Hom(E,F)=H\ar[d]^{\pi}&&\\
0\ar[r]&\Hom(C,C)\ar[r]^i& \End\CC^s\ar[ur]^{l_\mu}\ar[r]&
\CC^t\otimes \Ext^1(C,F)\ar[d] && \\
&&&0&&\\
}
\]
where $l_\mu(A)=\mu A$ and $r_\mu(B)=B\mu$.
Notice that the tangent space to the stabilizer of $\mu$ is
\[\textrm{T}(\Stab(\mu))=\{(A,B)\in\End \CC^s\times \End \CC^t | l_\mu(A)=r_\mu(B)\}.\]
We want to prove that $\dim\Stab(\mu)=\dim \textrm{T}(\Stab(\mu))=\dim\Hom(C,C)$.
Let us suppose that $A\in\End\CC^s$ satisfies
$l_\mu(A)\in \im(r_\mu)$.
Since the map $r_\mu$ is injective, there exists a unique $B\in\End
\CC^t$ such that
$(A,B)$ is in the stabilizer.
Moreover $\pi(l_\mu(A))=0$,
and thus, since the diagram is commutative, there
exists $\phi=i^{-1}(A)\in\Hom(C, C)$ which is unique, since $i$ is injective.
Conversely, we associate to every $\phi\in\Hom(C, C)$ a unique
$A=i(\phi)$. Since the sequences are exact and the diagram commutes,
we have $l_\mu(A)\in \Ker \pi=\im r_\mu$, i.e.\ there exists $B$ such that
the pair $(A,B)$ is in the stabilizer. Moreover, $B$ is unique, since
$r_\mu$ is injective by hypothesis $\Hom(C,F)=0$.

Hence, since this correspondence is one-to-one and linear,
it follows that $\dim \Stab(\mu)=\dim \Hom(C,C)$.
\end{proof}

\begin{theorem}\label{simpli}
Let $C$ be a generic bundle with resolution
\eqref{SD}. Then the following statements are equivalent:
\begin{enumerate}
\item[\rm{(i)}]
$C$ is simple,
\item[\rm{(ii)}]
$s^2-wst+t^2\leq 1$.
\end{enumerate}
\end{theorem}
\begin{proof}
From Theorem \ref{megakac} and Lemma \ref{lemmaStab} it follows that
$\rm{(ii)}$ implies  $ \rm{(i)}$.

To prove the other implication suppose that $C$ is simple.
Then, since $F$ is simple and $\Hom(F,C)\neq0$, it
follows that $\Hom(C,F)=0$. Hence applying Lemma
\ref{dimensioni}, we get that $\dim\Hom(C,C)=\dim\Stab(\mu)$.
Clearly
\[\dim\Stab(\mu)\geq\dim(\GL(s)\times\GL(t))-\dim H=s^2+t^2-wst,\]
thus
$1 =\dim\Hom(C,C)=\dim\Stab(\mu)\geq s^2+t^2-wst.$ Hence $\rm{(i)}$ implies
$\rm{(ii)}$.
\end{proof}

The following result will be used in Section \ref{stabilita}.

\begin{proposition}\label{deffo}
Assume that $E$ and $F$ are rigid,
and $\Ext^2(F,E)=0$.
If $s^2+t^2-wst\leq1$,
then the property of admitting resolution of the form
\eqref{SD} is invariant under small deformations.
\end{proposition}
\begin{proof}
Since $C$ is simple, the dimension of the space of matrices
in $\Hom(E^s,F^t)$ up to the action of $\GL(s)\times\GL(t)$ is
\[\dim H-\dim(\GL(s)\times\GL(t))+\dim\Stab(\mu)=wst-s^2-t^2+1\geq0.\]
We know that $\dim\Ext^1(C, C)\geq wst-s^2-t^2+1\geq0$.
On the other hand
we prove that $\dim\Ext^1(C,C)\leq wst-s^2-t^2+1$.
Indeed, by applying the functor $\Hom(C,-)$ to the resolution of $C$ we obtain
\begin{eqnarray*}
&0\rightarrow \Hom(C,E^s)\rightarrow \Hom(C,F^t)\rightarrow \Hom(C,C)\rightarrow \Ext^1(C,
E^s)\rightarrow &\\
&\rightarrow  \Ext^1(C,F^t)\rightarrow \Ext^1(C,C)\rightarrow \Ext^2(C, E^s)\rightarrow Q\rightarrow0.&
\end{eqnarray*}
Hence by the assumptions on $E$ and $F$ and the simplicity of $C$ it follows
\begin{eqnarray*}
&\dim\Ext^1(C,C)\leq s\left(\dim\Hom(C,E)-\dim\Ext^1(C,E)+
 \dim\Ext^2(C,E)\right)+&\\
& +t\left(\dim\Ext^1(C,E)-\dim\Hom(C,E)\right)+\dim\Hom(C,C)=-s^2+t(ws-t)+1,&
\end{eqnarray*}
which completes the proof.
\end{proof}

%
\section{Fibonacci bundles on $\PP^N$}
In this section we introduce the family of
{\em Fibonacci bundles}, which will replace exceptional bundles
in the canonical decomposition (see Section \ref{nonzimble}).
In fact these bundles satisfy some properties of
exceptional bundles, but in general they are not exceptional.

Given a bundle $G$, we denote by $G_z$ the fiber of $G$ at the point
$z\in\PP^N$.
Given a map of bundles $f:G \to L$, we denote by $f_z$ the restriction of the map
$f$ to the fiber at $z$, i.e.\  $f_z:G_z\to L_z$.
\begin{theorem}\label{mutazioni}
For any pair $(E,F)$ satisfying the basic hypotheses \eqref{ipo1} and \eqref{ipo2},
there exist the following sequences of bundles:
\begin{itemize}
\item
$0\rightarrow C_{n-1}\xrightarrow{i_{n}} C_{n}\otimes
W^{*}\xrightarrow{p_n} C_{n+1}\rightarrow 0\ $ if $n$ is odd,
\item
$0\rightarrow C_{n-1}\xrightarrow{i_{n}} C_{n}\otimes
W\xrightarrow{p_n} C_{n+1}\rightarrow 0\ $ if $n$ is even,
\end{itemize}
where
$C_0=E,\ C_1=F,$
and the map $i_n$ is recursively defined as follows:
\begin{itemize}
\item
$i_1:  E\to F\otimes W^{*}=F\otimes\Hom(E,F)^{*}$
is the canonical map and
\item
$i_{n}=(p_{n-1}\otimes id)\circ(id\otimes d): C_{n-1}\otimes
\CC\to C_{n-1}\otimes W\otimes W^{*}\to
C_{n}\otimes W$
(resp. $C_{n}\otimes W^{*}$),
where $d:\CC\to W\otimes W^{*}$ is the diagonal map.
\end{itemize}
We call the bundles $C_{n}$ ``Fibonacci bundles corresponding to $(E,F)$''.
\end{theorem}
\begin{proof}
In order to prove the theorem, we will go through the following
recursive steps for any $n$:
\begin{itemize}
\item[I.] we define the map
$i_n:C_{n-1}\to C_n\otimes W^{*}$ if $n$ is odd,
$i_n: C_{n-1}\to C_n\otimes W$ if $n$ is even;
\item[II.]
we prove by induction the property
\begin{eqnarray*}
(P_n)& \textrm{ for any } z\in \PP^N, \textrm{ for any }0\neq c\in
C_{n-1,z} \textrm{ the rank of }&\\
&i_{n,z}(c)\in\Hom(C_{n,z}^{*},W), \textrm{ resp. }\Hom(C_{n,z}^{*},W^{*}),
\textrm{ is bigger than }1;&
\end{eqnarray*}
\item[III.]
we prove that
$i_n$ is injective, i.e. that the rank of $i_n$ is constant;
\item[IV.]
we define $C_{n+1}:=\coker(i_n)$.
\end{itemize}

If $n=1$ the map $i_1$ is canonical, hence the property $(P_1)$ holds and
the fact that $E^*\otimes F$ is globally generated implies the
injectivity of $i_1$.

Now, let us assume the bundles $C_k$ to be defined for all $k\leq
n+1$, the map $i_k$ to be defined for all $k\leq n$, to satisfy $(P_k)$
and to be injective.

Let $n$ be odd. First, we define the map $i_{n+1}$.
By induction we have
\[0\rightarrow C_{n-1}\xrightarrow{i_n} C_n\otimes W^{*}\xrightarrow{p_n} C_{n+1}\rightarrow
0\]
where $p_n$ is
the projection induced by $i_n$.

By tensoring by $W$, we get the diagram
\[\xymatrix{
&&0\ar[d]&&\\
&&C_n\otimes \CC\ar[d]^-{id\otimes d}\ar@{=}[r]&
C_n\otimes \CC\ar[d]&\\
0\ar[r] & C_{n-1}\otimes W\ar@{=}[d]\ar[r]^-{i_n\otimes id}&
C_n\otimes W^{*}\otimes W\ar[r]^-{p_n\otimes id}\ar[d]
&C_{n+1}\otimes W\ar[r]&0\\
&C_{n-1}\otimes W\ar[r]
&C_n\otimes\Ad(W)\ar[d]^{i_n}&&\\
&&0&&}
\]
where $d:\CC\rightarrow W\otimes W^{*}$ is the diagonal map;
more explicitly if  $\{e_1,\ldots,e_w\}$ is a basis of $W$ and
$\{e_1^{*},\ldots,e_w^{*}\}$ the dual basis, then
\[d(1)=\sum_{i=1}^w e_i\otimes e_i^{*}.\]
We define the map $i_{n+1}$ as the following composition
\[i_{n+1}=(p_n\otimes id)\circ(id\otimes d): C_{n}\otimes
\CC\rightarrow C_{n+1}\otimes W.\]

Now we prove the property $(P_{n+1})$.
For any $z\in\PP^N$ we have $C_{n+1,z}=\frac{C_{n,z}\otimes
  W^{*}}{i_{n,z}(C_{n-1,z})}.$ Hence for any $c\in C_{n,z}$ we get
\[c\xrightarrow{(id\otimes d)_z} \sum_{i=1}^w c\otimes e_i \otimes
e_i^{*} \xrightarrow{(p_n\otimes id)_z} \sum_{i=1}^w\frac{c\otimes
  e_i^{*}}{i_{n,z}(C_{n-1,z})}\otimes e_i,\]
hence
\[i_{n+1,z}(c)=\sum_{i=1}^w\frac{c\otimes e_i^{*}}{i_{n,z}(C_{n-1,z})}\otimes
e_i.\]
If there exists $0\neq c\in C_{n,z}$ such that the rank of
$i_{n+1,z}(c)\in\Hom(C_{n+1,z}^{*},W)$ is $1$,
then for any $i\neq j$ there exist $\alpha_{ij},\beta_{ij}\in\CC$ such that
\[\alpha_{i,j}\frac{c\otimes e_i^{*}}{i_{n,z}(C_{n-1,z})}=\beta_{ij}\frac{c\otimes
  e_j^{*}}{i_{n,z}(C_{n-1,z})},\]
that is $\alpha_{ij}c\otimes e_i^{*}-\beta_{ij} c\otimes e_j^{*}= c\otimes
  (\alpha_{ij} e_i^{*}-\beta_{ij} e_j^{*})\in
i_{n,z}(C_{n-1})$, which contradicts $(P_n)$.
Therefore $(P_{n+1})$ is true.

Now in order to prove the injectivity of $i_{n+1}$, we show that
\begin{equation}\label{disgiu}
\im(i_n\otimes id)_z\cap \im(id\otimes d)_z=\{0\}, \textrm{ for all }z\in\PP^N.
\end{equation}
Indeed, for any $z\in\PP^N$, an element of $\im(id\otimes d)_z$ is of the form
$\sum_{i=1}^w c\otimes e_i \otimes e_i^{*}$ for some $c\in C_{n,z}$ and
if  $\sum_{i=1}^w c\otimes e_i \otimes e_i^{*} \in \im(i_n\otimes
id)_z$ then there exists an element
\[b=\sum_{i,j} \gamma_{ij} b_i\otimes e_j\in C_{n-1,z}\otimes W,\]
where $\{b_i\}$ is a basis of $C_{n-1,z}$ and
$\gamma_{ij}\in\CC$, such that
\[\sum_{i=1}^w c\otimes e_i^{*} \otimes e_i=(i_n\otimes id)_z(b).\]
It follows
\[\sum_{i=1}^w c\otimes e_i^{*}
\otimes e_i=\sum_{i,j} \gamma_{ij} i_{n,z}(b_i)\otimes e_j\]
and, projecting this equation on $e_j$, we get
\[c\otimes e_j^{*}=\sum_i \gamma_{ij} i_{n,z}(b_i)=i_{n,z}(\sum_i \gamma_{ij} b_i)\]
which contradicts $(P_n)$.
Hence \eqref{disgiu} is proved, and this implies the injectivity of
the map $i_{n+1}$ as a bundle map.

Finally we can define the bundle
$C_{n+2}:=\coker( i_{n+1}),$
and we get the exact sequence
\[0\rightarrow C_{n}\xrightarrow{i_{n+1}} C_{n+1}\otimes W\rightarrow C_{n+2}\rightarrow
0.\]

If $n$ is even, we repeat the same argument interchanging $W$ and
$W^{*}$ and this yields the following exact sequence
\[0\rightarrow C_{n}\xrightarrow{i_{n+1}} C_{n+1}\otimes W^{*}
\rightarrow C_{n+2}\rightarrow 0.\]
This completes the proof.
\end{proof}

\begin{theorem}\label{risoll}
For every $n\geq1$, a Fibonacci bundle $C_n$ on $\PP^N$ corresponding
to $(E,F)$ admits the following resolution
\begin{equation}\label{claim}
0\rightarrow E^{a_{n-1}}\rightarrow F^{a_n}\rightarrow
C_{n}\rightarrow 0,
\end{equation}
 with
$\{a_n\}=\{a_{w,n}\}$ as in \eqref{nufibo}.
\end{theorem}
\begin{proof}
We prove the statement by induction on $n$. If $n=1$ the
sequence \eqref{claim} is
$0\rightarrow F \rightarrow C_1\rightarrow 0,$ and the claim is true.

Now, we suppose that every $C_k$ admits a resolution of the form
\eqref{claim} for any $k\leq n$ and we  prove the same assertion
for $C_{n+1}$.
First we consider $n$ odd.
By the sequence
\[0\rightarrow C_{n-1}\rightarrow C_n\otimes W^{*}\rightarrow
C_{n+1} \rightarrow 0,\]
and by induction hypothesis, we have:
\[ \xymatrix{
       &      0          &  0                            &  &\\
0\ar[r]& C_{n-1}\ar[r]\ar[u]&   C_n\otimes W^{*}\ar[r]\ar[u]&   C_{n+1}\ar[r]&0 \\
       &F^{a_{n-1}}\ar[u]\ar[ru]^-{\alpha}   &   F^{a_n}\otimes W^{*}  \ar[u]&      & \\
       &E^{a_{n-2}}\ar[u]&   E^{a_{n-1}}\otimes W^{*}  \ar[u]&  & \\
       &    0\ar[u]            &  0    \ar[u]            &  &
} \]
where
we define the map $\alpha$
as  the composition of the known maps.

Since $\Ext^1(F,E)=0$, the map $\alpha$
induces a map $\widetilde\alpha:F^{a_{n-1}}\to F^{a_n}\otimes
W^{*}$ such that the diagram commutes.
Moreover if $\widetilde\beta$ is the restriction of $\widetilde\alpha$ to
$E^{a_{n-2}}$ the following diagram commutes:
\[ \xymatrix{
       &      0          &  0                            &  &\\
0\ar[r]& C_{n-1}\ar[r]\ar[u]&   C_n\otimes W^{*}\ar[r]\ar[u]&   C_{n+1}\ar[r]&0 \\
       &F^{a_{n-1}}\ar[u]\ar[ru]^-{\alpha}\ar[r]^-{\widetilde\alpha}   &   F^{a_n}\otimes W^{*}  \ar[u]&      & \\
       &E^{a_{n-2}}\ar[u]\ar[r]^-{\widetilde\beta}&   E^{a_{n-1}}\otimes W^{*}  \ar[u]&  & \\
       &    0\ar[u]            &  0    \ar[u]            &  &
} \]

This diagram implies that
$\Ker(\widetilde\alpha)\cong\Ker(\widetilde\beta)$, but since $E$ and
$F$ are simple, $\Ker(\widetilde\alpha)\cong F^a$ and
$\Ker(\widetilde\beta)\cong E^b$ for some $a,b\in\NN$.
Hence since $E$ and $F$ are indecomposable and $E\not\cong F$, by the
Krull-Schmidt theorem for vector bundles (see \cite{atiyah}),
we get
$\Ker(\widetilde\alpha)\cong\Ker(\widetilde\beta)=0.$
Thus
$\widetilde\alpha$ is injective and we can complete the diagram as follows:

\[ \xymatrix{
       &      0          &  0                            & 0 &\\
0\ar[r]& C_{n-1}\ar[r]\ar[u]&   C_n\otimes W^{*}\ar[r]\ar[u]&   \ar[u]C_{n+1}\ar[r]&0 \\
0\ar[r]&F^{a_{n-1}}\ar[u]\ar[ru]^-{\alpha}\ar[r]^-{\widetilde\alpha}
       &   F^{a_n}\otimes W^{*}  \ar[u]\ar[r]& F^{a_{n+1}} \ar[u]   \ar[r] &0 \\
0\ar[r]&E^{a_{n-2}}\ar[u]\ar[r]^-{\widetilde\beta}&
       E^{a_{n-1}}\otimes W^{*}  \ar[u]\ar[r]& E^{a_n} \ar[u]\ar[r]&0 \\
       &    0\ar[u]            &  0    \ar[u]            &0\ar[u]  &} \]
It follows that $C_{n+1}$ has the resolution
\[0\rightarrow E^{a_n}\rightarrow F^{a_{n+1}}\rightarrow C_{n+1}\rightarrow 0.\]

If we consider $n$ even, we replace $W$ with $W^{*}$
and we obtain the same result.
\end{proof}

We remark that it is possible to describe more explicitly the resolutions
of Fibonacci bundles.
Indeed for every $n\geq0$, a Fibonacci bundle $C_n$
corresponding to $(E,F)$ on $\PP^N$ has the following resolution
\[
0\rightarrow E\otimes A_n\rightarrow
F\otimes B_n\rightarrow C_{n}\rightarrow 0,
\]
where
\[A_{1}=0,\quad  B_1=\CC,\quad A_2=\CC,\quad B_2=W^{*},\]
and
\[
A_{n+1}=\frac{A_n\otimes W}{j_n(A_{n-1})} \textrm{ if } n \textrm{ even;}\quad\quad
 A_{n+1}=\frac{A_n\otimes W^*}{j_n(A_{n-1})} \textrm{ if } n \textrm{ odd;}
\]
\[
B_{n+1}=\frac{B_n\otimes W}{u_n(B_{n-1})} \textrm{ if } n \textrm{ even;}\quad\quad
B_{n+1}=\frac{B_n\otimes W^*}{u_n(B_{n-1})} \textrm{ if } n \textrm{ odd,}
\]
where $j_n$ and $u_n$ are recursively defined, in a similar way to the definition of $i_n$
in the statement of Theorem \ref{mutazioni}.

More explicitly, we define $j_1:0\to W$ as the zero map,
$u_1=d:\CC\to W^*\otimes W$ as the diagonal map.
For any $n\geq1$, we define $q_n$ and $r_n$ such that
\[0\rightarrow A_{n-1}\xrightarrow{j_{n}} A_{n}\otimes
U_n\xrightarrow{q_n} A_{n+1}\rightarrow 0 \quad\textrm{ and }\quad
0\rightarrow B_{n-1}\xrightarrow{u_{n}} B_{n}\otimes
U_n\xrightarrow{r_n} B_{n+1}\rightarrow 0,\]
where for brevity we denote $U_n=W$ if $n$ is even, $U_n=W^*$ if $n$ is odd.
Now we define, for any $n\geq2$,
\[j_{n}=(q_{n-1}\otimes id)\circ(id\otimes d): A_{n-1}\otimes
\CC\to A_{n-1}\otimes U_n\otimes U_n^*\to A_{n}\otimes U_n\]
and
\[u_{n}=(r_{n-1}\otimes id)\circ(id\otimes d): B_{n-1}\otimes
\CC\to B_{n-1}\otimes U_n\otimes U_n^*\to
B_{n}\otimes U_n.\]
\begin{remark}
It is easy to check that $A_k\cong B_{k-1}^*$ as
$\SL(V)$-representations, since all sequences of $\SL(V)$-modules
split. However it is possible that this isomorphism is not canonical,
because when $A_k$ and $B_k$ are decomposed as sums of irreducible
representations, some summand can appear with multiplicity bigger than
one.

In order to clarify the situation look at an example. Let $W=S^2V$
and $N=4$. We denote by
$\Gamma(a_1a_2a_3a_4)$ the irreducible representation of $\SL(V)$
with highest weight $\sum_ia_i\omega_i$, where $\omega_i$ are the
fundamental weights. With this notation we have for example
$A_3=S^2V=\Gamma(2000)$ and $B_2=S^2V^*=\Gamma(0002)$. Going
on, we can compute
\[ A_5=2\Gamma(3001) +2\Gamma(1101) +\Gamma(2000)+
\Gamma(0100)+\Gamma(4002)+\Gamma(2102) +\Gamma(0202)  \] and
\[
B_4=2\Gamma(1003) +2\Gamma(1011) +\Gamma(0002) +\Gamma(0010)+
\Gamma(2004) +\Gamma(2012) +\Gamma(2020).
\]
In this example it is evident that $A_5\cong
B_{4}^*$, nevertheless the isomorphism need not be canonic,
since two terms in the sum have multiplicity two.
\end{remark}


\begin{remark}
Since $a_{k-1}^2+a_k^2-wa_{k-1}a_k=1$, from Theorem \ref{simpli}
it follows that any generic bundle with resolution \eqref{claim} is simple.
In general this does not imply that any Fibonacci bundle is simple.
However with more assumptions on $(E,F)$, we will prove the simplicity of any Fibonacci bundle.
\end{remark}

\begin{lemma}\label{lemmuccio}
Assume that $E$ is rigid.
If $C_n$ is a
Fibonacci bundle corresponding to $(E,F)$, then
\[\dim\Hom(F,C_n)=a_n\quad\textrm{ and }\quad\dim\Hom(E,C_n)=a_{n+1}.\]
\end{lemma}
\begin{proof}
From Theorem \ref{risoll} we know that $C_n$ has resolution
\eqref{claim}.
Now, by applying respectively the functors $\Hom(F,-)$ and $\Hom(E,-)$ to this sequence,
we easily obtain that
$\dim\Hom(F,C_n)=a_n$ and
$\dim\Hom(E,C_n)=wa_n-a_{n-1}=a_{n+1},$ as claimed.
\end{proof}

Recall that by hypotheses $(R)$ we mean the union of conditions  \eqref{ipo1}, \eqref{ipo2}
and \eqref{ipo3}.

\begin{corollary}\label{simrig}
Assume that the pair $(E,F)$ satisfies the conditions $(R)$.
Then if a corresponding Fibonacci bundle $C_n$
is simple, it is also rigid.
\end{corollary}
\begin{proof}
Assume that $C_n$ is simple, i.e.\ $\dim\Hom(C_n,C_n)=1$. Since
$\Ext^1(F,F)=0$ and $\Ext^2(F, E)=0$, it follows
$\Ext^1(F,C_n)=0$.
Then by applying $\Hom(-,C_n)$ to the resolution \eqref{claim} and by Lemma \ref{lemmuccio}
we get
\[\dim\Ext^1(C_n,C_n)=\dim\Hom(C_n, C_n)- a_n\dim\Hom(F,C_n)+
a_{n-1}\dim\Hom(E,C_n)=\]
\[=1-a_n^2+a_{n-1}a_{n+1}=0,\]
hence $C$ is rigid.
\end{proof}


\begin{lemma}
\label{uno}
Assume that $(E,F)$ satisfies $(R)$ and, for any $n\geq0$, let $C_n$ be the corresponding
Fibonacci bundle. Then the following properties
$(I_n)$, $(II_n)$ and $(III_n)$ are satisfied for any $n\geq 1$:
\begin{eqnarray*}
(I_n) &\quad& \Hom(C_n,C_n)\cong\CC,\\
(II_n)&\quad& \Hom(C_{n},C_{n-1})=0,\quad \Ext^1(C_{n},C_{n-1})=0,\\
(III_n)&\quad &\Hom(C_{n-1},C_{n})\cong W,\quad\textrm{
  if }n\textrm{ odd,}\\
&\quad &\Hom(C_{n-1},C_{n})\cong W^*,\quad\textrm{
  if }n\textrm{ even.}
\end{eqnarray*}
\end{lemma}
\begin{proof}
We prove the lemma by induction on $n$.

Recall that $C_0=E$ and $C_1=F$. Hence by hypotheses \eqref{ipo1} and \eqref{ipo2}
we know that $C_1$ is simple, $\Hom(C_0,C_1)=W$ and
$\Hom(C_1,C_0)=\Ext^1(C_1,C_0)=0$.

Now suppose that $(I_k)$, $(II_k)$ and $(III_k)$ hold for all $k\leq
n$. Let us prove $(I_{n+1})$, $(II_{n+1})$ and $(III_{n+1})$.
First, suppose $n$ even. By the definition of Fibonacci bundles we get
the sequence
\begin{equation}
\label{mutazione}
0\rightarrow C_{n-1}\rightarrow C_n\otimes W\rightarrow
C_{n+1} \rightarrow 0.
\end{equation}
By applying the functor $\Hom(C_n,-)$ we get
\[0\rightarrow \Hom(C_{n},C_{n-1})\rightarrow
 \Hom(C_n,C_n)\otimes W\rightarrow
  \Hom(C_{n},C_{n+1})\rightarrow\Ext^1(C_{n},C_{n-1})\]
and from $(I_n)$ and $(II_n)$ we get
$\Hom(C_n,C_{n+1})\cong W$.
On the other hand if $n$ is odd, we consider the sequence
\[0\rightarrow C_{n-1}\rightarrow C_n\otimes W^*\rightarrow
C_{n+1} \rightarrow 0,\]
and we obtain with the same argument $\Hom(C_n,C_{n+1})\cong W^*$,
hence $(III_{n+1})$ follows.

Now, suppose $n$ even and let us apply $\Hom(-,C_n)$ to the sequence
\eqref{mutazione}:
\begin{eqnarray*}
&0\rightarrow\Hom(C_{n+1},C_{n})\rightarrow \Hom(C_n,C_n)\otimes
W^*\xrightarrow{\alpha} \Hom(C_{n-1},C_{n})\rightarrow&\\
&\rightarrow \Ext^1(C_{n+1},C_{n})\rightarrow\Ext^1(C_n,C_n)\otimes
W.&
\end{eqnarray*}
Since $\Hom(C_n,C_n)\cong \CC$ and $\alpha$ is the canonical
isomorphism $\Hom(C_{n-1},C_n)\cong W^*$,
we get $\Hom(C_{n+1},C_{n})=0$. Moreover from $(I_{n})$ we know that
$C_n$ is simple and by Corollary \ref{simrig} it follows that $C_n$ is
rigid, i.e.\ $\Ext^1(C_n,C_n)=0$. It implies that
$\Ext^1(C_{n+1},C_{n})=0$, and the property $(II_{n+1})$ follows.

Now let us prove $(I_{n+1})$.
First we apply $\Hom(-,C_{n-1})$
to \eqref{mutazione}, and we have
\[ \Hom(C_{n},C_{n-1})\otimes W^*
{\rightarrow} \Hom(C_{n-1},C_{n-1}) \rightarrow
\Ext^1(C_{n+1},C_{n-1})\rightarrow \Ext^1(C_{n},C_{n-1})\otimes
W^*, \] then  $(II_n)$ and $(I_{n-1})$ imply that
$\Ext^1(C_{n+1},C_{n-1})\cong\CC$. Finally by applying \linebreak
$\Hom(C_{n+1},-)$ to \eqref{mutazione}, we get
\[  \Hom(C_{n+1},C_{n})\otimes W\rightarrow
 \Hom(C_{n+1},C_{n+1})
\rightarrow\Ext^1(C_{n+1},C_{n-1})\rightarrow
 \Ext^1(C_{n+1},C_{n})\otimes W \]
and, using $(II_{n+1})$, we obtain
\[\Hom(C_{n+1},C_{n+1})\cong\Ext^1(C_{n+1},C_{n-1})\cong\CC,\]
hence $(I_{n+1})$ holds.
\end{proof}
As a consequence of Corollary \ref{simrig} and Lemma \ref{uno}, we obtain
the following theorem.

\begin{theorem}\label{gnulla}
For any bundle $C$ on $\PP^N$, with $N\geq2$, the following are equivalent:
\begin{enumerate}
\item[\rm{(i)}]
 $C$ is a Fibonacci bundle corresponding to some pair $(E,F)$
  satisfying hypotheses $(R)$,
\item[\rm{(ii)}] $C$ is simple and rigid.
\end{enumerate}
\end{theorem}
\begin{proof}
From Lemma \ref{uno} and from Corollary \ref{simrig} it follows that
property $\rm{(i)}$ implies $\rm{(ii)}$.
The other implication is easy to prove, because it suffices to choose the bundles $E=C(-d)$
and $F=C$, with $d\gg0$ such that the pair $(E,F)$ satisfies conditions $(R)$.
\end{proof}

\begin{remark}
Notice that, in particular, all the exceptional bundles are Fibonacci bundles with
respect to some pair $(E,F)$.
\end{remark}

\begin{lemma}\label{due}
Assume that $(E,F)$ satisfies $(R)$ and, for any $n\geq0$, let $C_n$ be the corresponding
Fibonacci bundle.
Then $\Ext^1(C_{n-1},C_n)=0$ for any $n\geq1$.
\end{lemma}
\begin{proof}
Let us apply $\Hom(-,C_n)$ to the resolution of $C_{n-1}$. Since $\Ext^1(F,C_n)=0$, by Lemma
\ref{lemmuccio}, we get
\begin{eqnarray*}
&\dim\Ext^1(C_{n-1},C_n)= \dim\Hom(C_{n-1},C_n)- a_{n-1}\dim\Hom(F,C_n)+&\\
&+a_{n-2}\dim\Hom(E,C_n)
=w-a_{n-1}a_n+a_{n-2}a_{n+1}=0.&
\end{eqnarray*}
\end{proof}
\section{Non-simple bundles}
\label{nonzimble}
In this section we investigate a generic bundle $C$ on $\PP^N$ with resolution \eqref{SD}
in the case $s^2+t^2-wst\geq 1$.
By Theorem \ref{simpli} we know that such a bundle $C$ is simple
only if $s^2+t^2-wst=1$, that is only if $C$  is a
deformation of a Fibonacci bundle.
Here we prove that when $s^2+t^2-wst\geq 1$ and the pair $(E,F)$ satisfies
hypotheses $(R)$, then any generic bundle $C$ is decomposable as a sum
of Fibonacci bundles.
In particular $C$ is simple if and only if it is a Fibonacci bundle (if
and only if $s^2+t^2-wst=1$).

\begin{remark}\label{ggg}
Since $E^*\otimes F$ is globally generated, we have
\[\rk E\rk F\leq w=\dim\HH^0(E^*\otimes F).\]
\end{remark}

The following lemma is a consequence of the second part of Theorem
\ref{megakac}. Here we give another elementary proof.
\begin{lemma}\label{bello}
For any
$s,t\in \NN$ satisfying $t\rk(F)-s\rk(E)\geq N$, and
\[s^2+t^2-wst\geq0,\]
there exist unique $k,n,m\in \NN$  such that the bundle
$C_k^{n}\oplus C_{k+1}^{m}$
admits a resolution of the form
\begin{equation} \label{res}
0\rightarrow E^s{\rightarrow} F^t\rightarrow C_k^{n}\oplus
C_{k+1}^{m}\rightarrow 0,
\end{equation}
 where $C_k$ and $C_{k+1}$ are
Fibonacci bundles.
\end{lemma}
\begin{proof}
By Remark \ref{ggg} and conditions $t\rk(F)-s\rk(E)\geq N$ and
$s^2+t^2-wst\geq0$, it follows that
$t\geq\big(\frac{w+\sqrt{w^2-4}}{2}\big)s$.

Fix $s,t$ such that ${t}>\big(\frac{w+\sqrt{w^2-4}}{2}\big)s$. Let
$\{a_k\}=\{a_{w,k}\}$ be the sequence
defined in \eqref{nufibo}. It is easy to check that the sequence
$\{\frac{a_{k+1}}{a_k}\}$ is decreasing to
$\frac{w+\sqrt{w^2-4}}{2}$. It follows that there exists $k\geq 1$
such that\[ \textrm{ either }\quad\frac{a_{k}}{a_{k-1}}=\frac{t}{s}\quad
\quad\textrm{ or }\quad\frac{a_{k+1}}{a_k}<\frac{t}{s}<\frac{a_k}{a_{k-1}}.\]
In the first case, since $(a_k,a_{k-1})=1$ by Remark \ref{rema},
there exists $n>1$ such that $t=na_k, s=na_{k-1}$, i.e.\ the
bundle $ C_k^{n}$ admits resolution \eqref{res}, with $m=0$. In
the second case, we solve the following system
\[\left\{
\begin{array}{l}
t=na_k+ma_{k+1},\\
s=na_{k-1}+ma_k.\\
\end{array}
\right.
\]
This system has discriminant $\Delta=a_k^2-a_{k+1}a_{k-1}=1$,
thus it admits a pair of integer solutions $(n,m)$. In
particular, $n>0$ because
 $\frac{t}{s}>\frac{a_{k+1}}{a_k}$, and
$m>0$ because $\frac{t}{s}<\frac{a_k}{a_{k-1}}$.
It follows that the bundle
$C_k^{n} \oplus C_{k+1}^{m}$ has resolution \eqref{res}.
\end{proof}
\begin{theorem}\label{decompo}
Let $(E,F)$ satisfy hypotheses $(R)$, and $s,t\in\NN$ satisfy
$t\rk(F)-s\rk(E)\geq N$.
Let $C$ be a generic bundle on $\PP^N$ with resolution \eqref{SD}.
Then
\[s^2+t^2-wst\geq0\quad\Rightarrow \quad C\cong C_k^{n}\oplus C_{k+1}^{m}\]
where $C_k$ and $C_{k+1}$ are Fibonacci bundles and $n,m\in\NN$ are unique.
\end{theorem}
\begin{proof}
It suffices to prove that the space of matrices $\mu\in H$ such that
$\coker(\mu)\cong C_k^{n}\oplus C_{k+1}^{m}$ is a dense subset of the
vector space $\Hom(E^s,F^t)$. Let us compute
$\dim\Ext^1(C_k^{n}\oplus C_{k+1}^m,C_k^{n}\oplus C_{k+1}^m)$.
By the property $(II)$ of Lemma \ref{uno}, by Theorem \ref{gnulla} and
Lemma \ref{due}, we obtain that $\dim\Ext^1(C_k^{n}\oplus
C_{k+1}^m,C_k^{n}\oplus C_{k+1}^m)=0$ for all $k\geq1$ and for all $n,m\in\NN$. Hence
the bundles $C_{k}^n\oplus C_{k+1}^m$ are rigid. It follows that the set of
matrices $\mu$ such that $\coker(\mu)$ is isomorphic to $C_k^n\oplus
C_{k+1}^m$ is open, hence dense in $H$. This completes the proof.
\end{proof}
\begin{remark}
Notice that a generic bundle which satisfies the hypotheses of
Theorem \ref{decompo} is rigid, hence homogeneous.
\end{remark}

\section{Some applications}
\subsection{First example}
If $E=\OO$ and $F=\OO(d)$ on $\PP^N$ (with $N\geq2$) the
bundles with resolution \eqref{primoesempio}
are cokernels of matrices of homogeneous polynomials of degree $d$.

Notice that hypotheses $(R)$ are satisfied
if either $N\geq3$, or $N=2$ and $1\leq d\leq 2$. Then all the
Fibonacci bundles on $\PP^N$ with either $N\geq3$ or $N=2$ and
$1\leq d\leq2$ are simple and rigid. Any generic deformation of a
Fibonacci bundle on $\PP^2$ with $d>2$ is simple. More precisely
we get the following classification:
\begin{corollary}
Let $C$ be a generic bundle on $\PP^N$
with resolution
\eqref{primoesempio}.
If either $N\geq3$ or $N=2$ and $1\leq d\leq2$, then
\[s^2+t^2-\binom{N+d}{d}st=1 \Leftrightarrow C \textrm{ is a Fibonacci bundle }
\Leftrightarrow C \textrm{ is simple and rigid}.\]
If $N=2$ and $d>2$,
\[s^2+t^2-\binom{N+d}{d}st=1 \Leftrightarrow C \textrm{ is a deformation of a
  Fibonacci}\Rightarrow C \textrm{ is simple}.\]
\end{corollary}
\begin{proof}
It is easy to check that if either $N\geq3$ or $N=2$ and $1\leq
d\leq2$, then $\dim\Hom(C,C)-\dim\Ext^1(C,C)=s^2+t^2-\binom{N+d}{d}st.$
\end{proof}

\begin{remark}
Let $C_k$ be a Fibonacci bundle on $\PP^N$ corresponding to $\OO,\OO(d)$
with $d\geq1$. Then
by applying the functor $\Hom(-,C_k)$ to the resolution of $C_k$,
we easily check the following properties:
\begin{eqnarray*}
&\Hom(C_k,C_k)=\CC,&\\
&\Ext^i(C_k,C_k)=0,\quad \textrm{ for all }\quad 1\leq i\leq
N-2,&\\
&\Ext^{N-1}(C_k,C_k)\cong \HH^{N}(\OO(-d))^{a_ka_{k-1}}.
\end{eqnarray*}
Thus all the Fibonacci bundles with $1\leq d\leq N$ are exceptional. In
particular if $d=1$, they are exactly the exceptional Steiner bundles studied
in \cite{mio}. The Fibonacci bundles with $d>N$ are not exceptional.
\end{remark}

\begin{corollary}
All the exceptional bundles on $\PP^N$ with resolution
\eqref{primoesempio} are exactly
the Fibonacci bundles corresponding to $\OO,\OO(d)$ for $1\leq d\leq N$.
All the semi-exceptional bundles on $\PP^N$ with resolution
\eqref{primoesempio} are of the form $C^n$, where $C$ is a Fibonacci
bundle as above.
\end{corollary}

\subsection{Bundles on $\PP^1$}
In this paper we have always supposed $N\geq2$, because the case
 $N=1$, corresponding to  bundles on $\PP^1$, is nowadays
trivial as it was solved by Kronecker in \cite{kronecker}.

In this case Theorem \ref{simpli} does not hold, since the fact that
$\dim\Stab(\mu)=1$ does not imply the simplicity of $\coker(\mu)$.  In
fact, since any bundle $C$ on $\PP^1$ is decomposable as a sum of
line bundles, $C$ is simple if and only if $C$ has rank $1$ if and
only if $C$ is exceptional.

On the other hand, there exists a canonical decomposition for all
bundles on $\PP^1$ with resolution \eqref{primoesempio}
for any $d\geq1$. Let us prove that a generic bundle with
resolution \eqref{primoesempio} is isomorphic to
$\OO(a)^n\oplus\OO(a+1)^m$. If $\frac{dt}{t-s}$ is
integer, then we choose $a=\frac{dt}{t-s}$, $n=t-s$ and $m=0$. If
$\frac{dt}{t-s}$ is not integer then we choose the unique integer
$\frac{dt}{t-s}-1<a<\frac{dt}{t-s}$. Then, as in the proof of
Theorem \ref{bello}, we see that the system
\[\left\{\begin{array}{l}
n+m=t-s\\
na+m(a+1)=dt\\
\end{array}\right.\]
admits a pair of integer positive solutions $(n,m)$. It follows
that the bundle $C=\OO(a)^n\oplus\OO(a+1)^m$ has resolution
\eqref{primoesempio}. Since $\dim\Ext^1(C,C)=0$, a generic bundle with
resolution \eqref{primoesempio} is isomorphic to $C$ and this implies
that there is a canonical reduction for matrices of polynomials in
two variables.
\begin{remark}
It follows that when we restrict any generic bundle on $\PP^N$ with
resolution \eqref{primoesempio}
to a generic $\PP^1\subset\PP^N$, the splitting type is of the form
$\OO(a)^n\oplus\OO(a+1)^m$, hence
it is as balanced as possible. %
\end{remark}

\subsection{Second example}
Given $0<p<N$, let us consider one of the following
pairs of bundles on $\PP^N=\PP(V)$, with $N\geq2$,
\begin{itemize}
\item $E=\OO(-1)$ and $F=\Omega^p(p)$,
\item $E=\Omega^p(p)$ and $F=\OO$.
\end{itemize}
It is easily seen that in these two cases hypotheses $(R)$ hold. Then
we can apply Theorems \ref{simpli} and \ref{decompo} to the corresponding
cokernel bundles and
we get the following consequences.

\begin{corollary}
Given $0<p<N$, let $C$ be a generic bundle with resolution either \eqref{secondoes}
or \eqref{secondoes-bis}.
Then
\begin{itemize}
\item
$C$ is simple $\quad \Leftrightarrow\quad s^2-wst+t^2\leq 1$,
\item
$C$ is simple and rigid $\quad\Leftrightarrow\quad$ C is a Fibonacci bundle $C_k$,
\item
$s^2-wst+t^2\geq 1\quad \Rightarrow\quad
C\cong C_k^n\oplus C_{k+1}^m$ for unique $n,m,k\in \NN$,
\end{itemize}
where $w=\binom{N+1}{N-p}$ in case \eqref{secondoes}, and
$w=\binom{N+1}{p}$ in case \eqref{secondoes-bis}.
\end{corollary}

Notice that also
exceptional Steiner bundles belong to this class. More precisely we have:
\begin{proposition}
Any exceptional Steiner bundle $S_{k+1}$ on $\PP^N$ of the form 
\[0\rightarrow \OO(-2)^{a_{k}}\rightarrow \OO(-1)^{a_{k+1}}\rightarrow
S_{k+1} \rightarrow 0\]
is isomorphic to a Fibonacci bundle $C_k$ associated to
the pair  $\OO(-1)$, $\Omega^{N-1}(N-1)$, i.e.\
\[0\rightarrow \OO(-1)^{a_{k-1}}\rightarrow
\Omega^{N-1}(N-1)^{a_k}\rightarrow C_k \rightarrow 0.\]
\end{proposition}
\begin{proof}
It suffices to apply the theorem of Be{\u\i}linson (see for example
\cite{AOforum}) to the bundle $S_{k+1}$.

Let us compute the dimension $\hh^i(F(-j))$ for
$0\leq i,j\leq N$.
From the resolution of $S_{k+1}$, it is easily seen
that $\hh^i(S_{k+1}(-j))=0$ for any $i\geq0$ and $0\leq j\leq N-2$.
Moreover $\hh^i(S_{k+1}(-N+1))=0$ for any $i\neq N-1$, and
$\hh^{N-1}(S_{k+1}(-N+1))=a_k.$
Finally $\hh^i(S_{k+1}(-N))=0$  for any $i\leq N-2$, and
$\hh^{N-1}(S_{k+1}(-N))-\hh^{N}(S_{k+1}(-N))=a_{k-1}$.
By Serre duality we know that $\hh^{N}(S_{k+1}(-N))=\hh^0(S_{k+1}^*(N-N-1))$.
Since $S_{k+1}$ is simple and $\HH^0(S_{k+1}(1))\neq0$, then
$\HH^0(S_{k+1}^*(-1))=0$. Hence $\hh^{N-1}(S_{k+1}(-N))=a_{k-1}$.

Then by applying the theorem of Be{\u\i}linson,
we get that $S_{k+1}$ admits the resolution
\[0\rightarrow \OO(-1)^{a_{k-1}}\rightarrow
\Omega^{N-1}(N-1)^{a_k}\rightarrow S_{k+1} \rightarrow 0.\]
Hence by the rigidity of the Fibonacci bundles we conclude that $S_{k+1}\cong C_k$.
\end{proof}

\subsection{Third example}
On $\PP^2$, we consider $E=S^pQ$ and $F=S^rQ(d)$ and the bundles
with resolution \eqref{terzoesempio}
for $d>p+1$.
Hypotheses $(R)$ are never satisfied, since $S^pQ$ is not rigid.
Hence we can deduce only the following result.
\begin{corollary}
Given $r,p\geq 1$ and $d>p+1$, a generic bundle $C$ with resolution
\eqref{terzoesempio}
is simple if and only if $s^2-wst+t^2\leq 1$, where
$w=\dim\HH^0(S^pQ\otimes S^rQ\otimes\OO(d-p))$.
\end{corollary}


\section{Stability}\label{stabilita}
In this last section we present some results about stability: first
we consider the
exceptional Steiner bundles on $\PP^N$ with $N\geq2$, then we
restrict our attention to
 bundles on $\PP^2$ and we utilize some important
results of Dr\'ezet and Le Potier.

\subsection{Stability of exceptional Steiner bundles on $\PP^N$}
Here we prove the stability of the Steiner
exceptional bundles on $\PP^N$ for any $N\geq2$.
Recall that on $\PP^N$, with $N>3$, the general problem of the stability of
exceptional bundles is still open.

\begin{theorem}\label{stabbi}
Any exceptional Steiner bundle $S_n$ on $\PP^N$ is stable for all
$n\geq0$ and for any $N\geq2$.
\end{theorem}
\begin{proof}
We prove the theorem by induction on $n$. If $n=0,1$, we get
$S_0=\OO$ and $S_1=\OO(1)$, which are stable, since they are line
bundles. Let us suppose that $S_k$ is stable for all $k\leq n$ and
let us prove the stability of $S_{n+1}$.
Assume by contradiction that $S_{n+1}$ is not stable. Then
there exists a quotient $Q$ such that
\[\mu(Q)\leq\mu(S_{n+1}).\]
We can suppose that $Q$ is stable. From Theorem \ref{mutazioni} we know that there exists the sequence
\[0\rightarrow S_{n-1}\rightarrow S_n\otimes U_n\rightarrow
S_{n+1} \rightarrow 0,\] where $U_n=W$ if $n$ is even, $U_n=W^*$
if $n$ is odd. It follows that $Q$ is also a quotient of
$S_n\otimes U_n$ and so, from the stability of $S_n$, we obtain
\[\mu(Q)\geq\mu(S_{n}).\] From the resolution of exceptional
Steiner bundles
\[0\rightarrow  \OO^{a_{n-1}}\rightarrow
\OO(1)^{a_n}\rightarrow S_{n}\rightarrow 0,\] we compute
$\mu(S_n)=\frac{a_n}{a_n-a_{n-1}}$. it is easy to check that
\[\mu(S_n)=\frac{a_n}{a_n-a_{n-1}}<\frac{a_{n+1}}{a_{n+1}-a_{n}}=\mu(S_{n+1})\]
and, denoting $r_k=a_k-a_{k-1}$, we compute
\[\frac{a_{n+1}}{r_{n+1}}-\frac{a_n}{r_n}=
\frac{1}{r_{n+1}r_n}.\]
 Hence, denoting by $\frac{c}{r}$ the slope of $Q$, we have to find
 two positive
 integer $c,r$ such that $r<r_{n+1}$ and
\[\frac{a_n}{r_n}\leq\frac{c}{r}\leq
\frac{a_n}{r_n}+\frac{1}{r_{n+1}r_n}.\]
 With simple computations we get
\[0\leq\frac{r_nc-a_nr}{r}\leq
\frac{1}{r_{n+1}}\] and, since $r<r_{n+1}$, the only possibility
is ${r_nc-a_nr}=0$, i.e.\ $\mu(Q)=\mu(S_n)$.
Now since $S_n\otimes U_n$ is polystable (in fact it is the direct sum
of $N+1$ copies of the stable bundle $S_n$) and $Q$ is stable with the same slope,
it follows that $Q=S_n$. Then $S_n$ has to be a quotient of $S_{n+1}$
and this is impossible because $\Hom(S_{n+1},S_n)=0$. This completes the proof.
\end{proof}

\subsection{Stability of bundles on $\PP^2$}
The problem of the stability of vector bundles on $\PP^2$
has been studied by Dr\'ezet and Le Potier.
In particular, in
\cite{DrezetLePotier} they found a criterion to check the
existence of a stable bundle with given rank and Chern classes,
but this criterion is complicated to apply even for Steiner
bundles.

Moreover, from another result of Dr\'ezet (see Theorem $3.1$ of
\cite{Drezet99}) we know that if there exist no semi-stable bundles
with given rank and Chern classes, then the generic bundle in the
space of {\em prioritary bundles} with these rank and Chern classes is
decomposable, hence non-simple.

A vector bundle $P$ on $\PP^2$ (or a coherent torsionfree sheaf) is called
{\em prioritary} when \[\Ext^2(P,P(-1))=0.\]
Prioritary bundles were introduced by Hirschowitz and Laszlo in
\cite{HirschowitzLaszlo}.

It is easily seen that if $E$ and $F$ are prioritary and $\Ext^1(E,
F(-1))=0$  then any bundle $C$ on $\PP^2$  with resolution
\eqref{SD} is prioritary.
On the other hand, if the pair $(E,F)$ satisfies hypotheses $(R)$,
by Proposition \ref{deffo} we get that a generic cokernel
bundle $C$ in \eqref{SD} is also generic in the space of prioritary bundles.

This implies that our Theorem \ref{simpli} provides a criterion
for the stability of generic bundles $C$ on $\PP^2$ with resolution
\eqref{SD},
for any $(E,F)$ satisfying hypotheses $(R)$.
Precisely we get the following result:
\begin{theorem}
Let $E$ and $F$ be two prioritary bundles on $\PP^2$ satisfying $(R)$ and such
 that $\Ext^1(E,F(-1))$. Let $C$ be defined by  \eqref{SD}.
 Then the following statements are equivalent:
\begin{enumerate}
\item[\rm{(i)}] $C$ is stable,
\item[\rm{(ii)}] $s^2-wst+t^2\leq 1$.
\end{enumerate}
\end{theorem}
\begin{remark}
In particular the previous theorem implies that any
Fibonacci bundle on $\PP^2$ with respect to $(E,F)$ is stable,
if $(E,F)$ satisfies $(R)$ and $\Ext^1(E,F(-1))$.
\end{remark}
\begin{remark}
If the pair $(E,F)$ satisfies only the basic hypotheses \eqref{ipo1} and \eqref{ipo2},
then we obtain the
stability of a generic deformation of a corresponding Fibonacci bundle
in the space of prioritary bundles.
\end{remark}

In the particular case $E=\OO$ and $F=\OO(d)$,
we have the following result:
\begin{theorem}\label{ciao}
Let $d=1,2$ and $C$ be a generic bundle on
$\PP^2$ defined by the exact sequence
\eqref{primoesempio}.
Then the following statements are equivalent:
\begin{enumerate}
\item[\rm{(i)}]$C$ is stable,
\item[\rm{(ii)}] $s^2-\binom{d+2}{2}st+t^2\leq 1$,
\item[\rm{(iii)}] either $C$ is a Fibonacci (and exceptional) bundle
or $t\leq \frac{1}{2}\left(\binom{d+2}{2}+\sqrt{\binom{d+2}{2}^2-4}\right)s$.
\end{enumerate}
\end{theorem}
We stress that this criterion is equivalent to the Dr\'ezet-Le
Potier criterion in the particular case of bundles with resolution
\eqref{primoesempio}. Nevertheless our proof is completely independent, and
it seems difficult to deduce it directly from \cite{DrezetLePotier}.

From the description of
non-simple bundles (Theorem \ref{decompo}), we can classify all the
strictly semi-stable bundles on $\PP^2$ with resolution \eqref{primoesempio}
and we get the following result.
\begin{corollary}
Let $F$ be a generic bundle on $\PP^2$ with resolution
\eqref{primoesempio} with $d=1,2$. Then $F$ is strictly semi-stable if and
only if it is the sum of $n>1$ copies of a Fibonacci bundle, if and
only if it is semi-exceptional.
\end{corollary}

Finally we remark that the results of this section allow us to improve a
theorem of Hein, contained in the appendix of \cite{Brenner-Hein}, about
the stability of a generic {\em syzygy bundle}, i.e. of a bundle $G$
on $\PP^N$ with resolution
\begin{equation}\label{syz}
0\rightarrow \OO\rightarrow \OO(d)^t\rightarrow G\rightarrow 0.
\end{equation}
In fact Theorem $A.1$ of \cite{Brenner-Hein} gives a sufficient
condition for the semi-stability of syzygy bundles on $\PP^N$ ($t\leq d(N+1)$) and
Theorem $A.2$ for the stability of syzygy bundles on $\PP^2$.
In particular he proves that a sufficient condition for the stability
of a syzygy bundle with resolution \eqref{syz}
is \[t\leq\frac{4}{5}d+1.\]
The following improvement is a consequence of our Theorem \ref{ciao}.
\begin{corollary}
A generic bundle $G$ with resolution \eqref{syz} on $\PP^2$ for $d=1,2$ is stable
if and only if $t\leq 3d$.

\end{corollary}

{\bf {Acknowledgement.}}
I would like to thank Giorgio Ottaviani, for his precious advice
and for many useful discussions.
%
%
I also thank the referee
for providing many suggestions that improved the exposition of the paper.
%
%
This work was supported by the MIUR in the framework of the National
Research Project ``Propriet\`a geometriche delle variet\`a reali
e complesse''.

%
%


\begin{thebibliography}{10}

\bibitem{AOforum} {Vincenzo Ancona and Giorgio Ottaviani, {Some
    applications of {B}e{\u\i}linson's
  theorem to projective spaces and quadrics}, Forum Math. \textbf{3} (1991),
  no.~2, 157--176}.

\bibitem{atiyah}
Michael Atiyah, {On the {K}rull-{S}chmidt theorem with application to
  sheaves}, Bull. Soc. Math. France \textbf{84} (1956), 307--317.

\bibitem{mio}
Maria~Chiara Brambilla, {Simplicity of generic {S}teiner bundles},
Bollettino Unione Matematica Italiana (8), \textbf{8-B} (2005), 723--735.

\bibitem{tesi}
Maria~Chiara Brambilla, {Simplicity of vector bundles on ${\PP}^{n}$ and exceptional
  bundles}, PhD thesis, University of Florence, 2004.

\bibitem{Brenner-Hein}
Holger Brenner, {Looking out for stable syzygy bundles}, arXiv:
  math.AG/0311333 (2005).

\bibitem{Drezet99}
Jean-Marc Dr{\'e}zet, {Vari\'et\'es de modules alternatives}, Ann. Inst.
  Fourier (Grenoble) \textbf{49} (1999), no.~1, v--vi, ix, 57--139.

\bibitem{DrezetLePotier}
Jean-Marc Dr{\'e}zet and Joseph Le~Potier, {Fibr\'es stables et fibr\'es
  exceptionnels sur {${\PP}^2$}}, Ann. Sci. \'Ecole Norm. Sup. (4) \textbf{18}
  (1985), no.~2, 193--243.

\bibitem{GoroRudakov}
Alexei~L. Gorodentsev and Alexei~N. Rudakov, {Exceptional vector bundles
  on projective spaces}, Duke Math. J. \textbf{54} (1987), no.~1, 115--130.

\bibitem{HirschowitzLaszlo}
Andr{\'e} Hirschowitz and Yves Laszlo, {Fibr\'es g\'en\'eriques sur le
  plan projectif}, Math. Ann. \textbf{297} (1993), no.~1, 85--102.

\bibitem{kac}
Victor~G. Kac, {Infinite root systems, representations of graphs and
  invariant theory}, Invent. Math. \textbf{56} (1980), no.~1, 57--92.

\bibitem{kronecker}
Leopold Kronecker, {Algebraische reduction der schaaren bilinearer
  formen}, 763--776, Berlin: S.-B. Akad., 1890.

\bibitem{Rudakovseminari}
Alexei~N. Rudakov, {Helices and vector bundles}, London Math. Soc. Lecture
  Note Series, vol. 148, Cambridge Univ. Press, Cambridge, 1990.

\bibitem{Zube}
Severinas~K. Zube, {The stability of exceptional bundles on
  three-dimensional projective space}, Helices and vector bundles, London Math.
  Soc. Lecture Note Series, vol. 148, Cambridge Univ. Press, Cambridge, 1990,
  pp.~115--117.

\end{thebibliography}
\end{document}